\documentclass{article}

\usepackage{amsmath,amsfonts,latexsym,theorem, 
}
\nonstopmode
\usepackage{hyperref}
\usepackage[T1]{fontenc}
\usepackage[utf8]{inputenc}
\usepackage{authblk}
\usepackage{amssymb}
\usepackage{mathtools}
\usepackage{indentfirst}
\usepackage{cite}
\usepackage{tikz}
\usetikzlibrary{matrix, shapes.arrows}
\usetikzlibrary{arrows.meta}

\numberwithin{equation}{section}

\newtheorem{Theorem}{Theorem}[section]

\newtheorem{Lemma}{Lemma}[section]

\newenvironment{Proofc}[1]{\smallskip\par\noindent\textsc{#1}\quad}%
{\hfill$\Box$\bigskip\par}
\newenvironment{Proof}{\begin{Proofc}{Proof}}{\end{Proofc}}

\title{Infinite combinatorial Ricci flow in spherical background geometry} 
\author{Chang Li, Yangxiang Lu and Hao Yu}

\date{}

\begin{document}
\maketitle

\begin{abstract}
Since the fundamental work of Chow-Luo \cite{CL03}, Ge \cite{Ge12,Ge17} et al., the combinatorial curvature flow methods became a basic technique in the study of circle pattern theory. In this paper, we investigate the combinatorial Ricci flow with prescribed total geodesic curvatures in spherical background geometry. For infinite cellular decompositions, we establish the existence of a solution to the flow equation for all time. Furthermore, under an additional condition, we prove that the solution converges as time tends to infinity. To the best of our knowledge, this is the first study of an infinite combinatorial curvature flow in spherical background geometry.
\end{abstract}

\section{Introduction}
\subsection{Background}
Circle patterns serve as discrete analogs of conformal structures on Riemann surfaces, a concept reintroduced by Thurston  \cite{st, Th76}. At the 1985 International Symposium celebrating the proof of the Bieberbach Conjecture, Thurston conjectured that  a sequence of functions, which are correspondences of circles within two constructed infinite hexagonal circle patterns, could be used to approximate the Riemann mapping from a bounded simply connected domain to the unit disk.  This conjecture was subsequently proven by Rodin-Sullivan\cite{ro2}.

The problem of finding a circle pattern on a surface is equivalent to determining a circle pattern metric with vanishing discrete Gaussian curvature. This area of research has been significantly advanced through the application of Ricci flow, a fundamental tool in geometric analysis originally developed by Hamilton \cite{H82}. The Ricci flow technique has proven instrumental in studying geometric structures on a given manifold, playing a pivotal role in differential geometry. Notably, Perelman \cite{Perelman} successfully resolved both the Poincaré conjecture and Thurston's geometrization conjecture by using the Ricci flow. Further advancements  on geometric structures related to the Ricci curvature have been achieved in recent work, including Jiang-Naber \cite{JN21} and Cheeger-Jiang-Naber \cite{CJN21}.

Chow-Luo \cite{CL03} introduced a discrete version known as the combinatorial Ricci flow through the evolution equation for finding circle pattern metrics with zero discrete Gaussian curvatures on triangulations of compact surfaces: 
\begin{equation}
	   \frac{\mathrm{d} r_i}{\mathrm{d} t} = -K_i\, s(r_i),
   \end{equation}  
where the function $s(r)$ is given by $s(r) = r$ in Euclidean background geometry and $s(r) = \sinh r$ in hyperbolic background geometry. Here, $K_i = 2\pi - \alpha_i$ represents the discrete Gaussian curvature at the vertex $i$, where $\alpha_i$ denotes the cone angle at the center of the circle $C_i$. They established the long-term existence and the convergence of the flow under some combinatorial conditions in both Euclidean and hyperbolic background geometries. Furthermore, their work on the combinatorial Ricci flow yielded an alternative proof of Thurston's famous circle pattern theorem. By proving that the solutions to the combinatorial Ricci flows converge exponentially fast to Thurston's circle patterns on surfaces, they established a deep connection between the combinatorial Ricci flow and circle patterns. This work introduced the combinatorial Ricci flow as an efficient framework for studying geometric structures on surfaces.

The combinatorial Ricci flow has found wide applications, notably in geometric topology (see \cite{GHZ21,GHZ21-2,GJ16,GJ17,GJ17-2,GJS22,GX16,GX16DGA,GX18}). In particular, the work of Feng-Ge-Hua \cite{FGH21}, and Ge-Hua \cite{GH20} had successfully applied the combinatorial flow methods to the study of 3-dimensional geometric topology, which are pioneering and fundamental. Furthermore, to address related problems, the combinatorial Ricci flow has been appropriately extended, such as the combinatorial Calabi flow for Thurston's circle pattern metrics (see \cite{Ge12}). Subsequent work in \cite{Ge17, GH18, GX16DGA} showed that these flows exist globally and converge exponentially fast to Thurston's circle patterns on surfaces if and only if there exists a circle pattern metric with constant (resp. zero) discrete Gaussian curvature in Euclidean (resp. hyperbolic) background geometry, establishing a basic framework and essential techniques in this field. 

Prior research has primarily focused on combinatorial curvature flows for polyhedral surfaces in Euclidean and hyperbolic background geometries. In contrast, results on circle patterns in spherical background geometry remain comparatively scarce. This disparity stems from the foundational work of Chow-Luo \cite{CL03}, which builds upon a functional introduced by Colin de Verdière \cite{Ve91}. While this functional exhibits convexity in both Euclidean and hyperbolic background geometries, it fails to maintain this property in spherical background geometry. Recently, Nie \cite{Nie} proposed a novel functional incorporating total geodesic curvatures that remains convex even in spherical background geometry. Utilizing this enhanced functional, Nie established both existence and rigidity results for circle patterns in spherical geometry with prescribed total geodesic curvatures at each vertex. Then, Ge-Hua-Zhou \cite{GHZ24} developed the first combinatorial curvature flow in spherical background geometry and established its refined properties. For convenience, we refer to this flow as the combinatorial Ricci flow in spherical background geometry. 
 
In this paper, we study spherical metrics (possibly with singularities) on general surfaces, including noncompact surfaces. Specifically, we focus on a constant curvature metric for a given surface S, where the metric may have singularities at certain points. As is well-known, there have been many results for the finite setting on circle pattern theory that are used to construct metrics (see \cite{BBI22}). In particular, infinite circle patterns provide more details for general surfaces, which is a direct generalization of finite circle patterns. In Euclidean and hyperbolic background geometries, there are particularly profound results in the infinite setting (see \cite{H82}). However, a significant gap remains in the study of infinite circle patterns in spherical background geometry. In this work, we show several properties of spherical metrics with infinite circle patterns for arbitrary surfaces S.
 
The study of circle patterns on infinite triangulations of noncompact surfaces constitutes an area of significant mathematical interest. Indeed, the infinite setting presents substantial challenges, as techniques developed for the finite setting cannot be directly generalized. The foundational rigidity result in this setting was established by Rodin-Sullivan \cite{ro2}. Subsequent developments include classification results for circle patterns on infinite triangulations of the complex plane $\mathbb{C}$ obtained by He-Schramm \cite{He95}, later extended by He \cite{He99}. Recently, Ge-Hua-Zhou \cite{GHZ25} studied the combinatorial Ricci flow on infinite disk triangulations. These fundamental contributions primarily employ techniques from the elliptic partial differential equation and conformal geometry. A natural and intriguing direction is to explore whether parabolic methods can also be effectively employed to study circle patterns on infinite triangulations. We note that the results in \cite{GHZ25,He95,He99} specifically address infinite circle patterns in either Euclidean or hyperbolic background geometry. To our knowledge, our work presents the first results for infinite circle patterns in spherical background geometry.

\subsection{The infinite circle patterns}
Consider a graph $G$ with an edge function $\Theta: E \rightarrow[0, \pi / 2]$. The following problem arises: Does there exist a circle pattern $\mathcal{P}$ whose contact graph is combinatorially equivalent to $G$ and whose intersection angles are given by $\Theta$? Furthermore, if such a pattern exists, what can be said about the uniqueness?

This problem is well-defined when $0 \leq \Theta(e) \leq \pi / 2$ for all edges $e \in E$. For finite graphs, the existence question is fully resolved by Thurston's interpretation of Andreev's theorem. The uniqueness, however, is only fully characterized when $G$ is the 1-skeleton of a triangulation of the 2-sphere (see \cite{MR89,Th80}).

We now state the Koebe uniformization of circle patterns. Consider a connected planar graph $G$ with no loops or multiple edges. A vertex subset $S \subseteq$ $V(G)$ is called a separating set if there exist at least two vertices in $V(G) \backslash S$ such that every path connecting them intersects $S$. Given an edge function $\Theta: E(G) \rightarrow[0, \pi / 2]$, by Thurston's interpretation of Andreev's theorem and a compactness argument, we know that the existence of a circle pattern which realizes the data $(G, \Theta)$  is equivalent to the following conditions when $G$ has at least five vertices (see proof in \cite{Th80}):\\\\
(C1) If a simple loop in $G$ formed by three edges $e_1, e_2, e_3$ separates the vertices of $G$, then $\sum_{i=1}^3 \Theta\left(e_i\right)<\pi$;\\\\
(C2) If $v_1, v_2, v_3, v_4=v_0$ are distinct vertices in $G$ and if $\{v_{i-1}, v_i\} \in E$ and $\Theta\left(\{v_{i-1}, v_i\}\right)=\pi / 2$, $i=1, \ldots, 4$, then either $\{v_0, v_2\}$ or $\{v_1, v_3\}$ is an edge in $G$.\\\\
Besides, it is well-known that the existence of a circle pattern is equivalent to the long-term existence and convergence of the combinatorial Ricci flow. As we previously mentioned, the theory of the combinatorial Ricci flow in the finite setting has been well established. However, the infinite setting requires stronger hypotheses beyond these basic conditions (C1) and (C2). The following fundamental result was established by He \cite{He99}:
\begin{Theorem}\label{Thm He}
Let $G$ be a disk triangulation graph, and let $\Theta: E \rightarrow[0, \pi / 2]$ be a function defined on the set of edges. Assume that conditions $\rm (C1)$ and $\rm (C2)$  hold.\\\\
(i) If $G$ is VEL-parabolic, then there is a locally finite disk pattern in $\mathbb{C}$ which realizes the data $(G, \Theta)$.\\\\
(ii) If $G$ is VEL-hyperbolic, then there is a locally finite disk pattern in $U$ which realizes $(G, \Theta)$.
\end{Theorem}

Hence, it is necessary to study the long-term existence and convergence of the flow in the infinite setting. However, in this setting, classical ODE theory becomes inapplicable, which is the fundamental difficulty. Our approach involves constructing an exhaustive sequence of finite domains and employing the associated solutions to approximate the solution to the flow. Moreover, it is worth emphasizing that the relationship between the existence and convergence of the combinatorial flow in the infinite setting and condition `VEL' presents a significant direction for future research.

\subsection{Main results}
Let $\mathcal{D}=(V, E, F)$ be an infinite and locally finite cellular decomposition, where $V, E$ and $F$ represent the sets of vertices, edges and faces, respectively. We define a circle pattern metric $r: V \rightarrow (0, \frac{\pi}{2})$, an intersection angle $\Theta: E \rightarrow(0, \frac{\pi}{2}]$ and total geodesic curvature $T: V \rightarrow \mathbb{R}$. Note that the total geodesic curvature $T_i=\alpha_i - c\operatorname{Area} (C_i)$ by Gauss-Bonnet formula, where $\alpha_i$ is the cone angle, $C_i$ is the circle centered at $v_i$ with radius $r_i$ and $c=1$ in spherical background geometry, $c=0$ in Euclidean background geometry, $c=-1$ in hyperbolic background geometry (see later sections for the formal definition). Since we use circle patterns to describe the metric structure of the  geometric models, compared to $\alpha_i$, total geodesic curvature can be seen as a modified cone angle singularity which describes another possible discrete analogue to continuous geometric model. The following conditions establish connections between geometry and combinatorics.\\\\
(S1)   $\hat {T}_v >0, \quad \forall v \in V;$\\\\
(S2)   $\ \sum_{v \in U} \hat{T}_v<2 \sum_{e \in E(U)} \Theta(e), \quad \forall U \subset V;$\\\\
(S3)   $T(r(0))\geq \hat{T}.$\\\\
Now, we present our main results. For the combinatorial Ricci flow equation \eqref{flowequation}  with prescribed total geodesic curvatures, we have the main theorems as follows,
\begin{Theorem}[Existence]\label{exist}
	Let $\mathcal{D}=(V, E, F)$ be an infinite cellular decomposition with intersection angle $\Theta \in(0, \frac{\pi}{2}]^E$. For any initial value $r(0)\in (0, \frac{\pi}{2})^V$, and $\{\hat{T}_i\}_{v_i \in V}$ that satisfy condition $\rm (S1)$ and $\rm (S2)$, there exists a solution $r(t)\in (0, \frac{\pi}{2})^V$ for $ t \in[0, \infty)$, to the flow \eqref{flowequation}.
\end{Theorem}

\begin{Theorem}[Convergence]\label{conv}
	Let $\mathcal{D}=(V, E, F)$ be an infinite cellular decomposition with intersection angle $\Theta \in(0, \frac{\pi}{2}]^E$. For an initial value $r(0)\in (0, \frac{\pi}{2})^V$, and $\{\hat{T}_i\}_{v_i \in V}$ that satisfy condition $\rm (S1)$, $\rm (S2)$ and $\rm (S3)$, there exists a solution $r(t)$ to the flow \eqref{flowequation}  such that $r(t)$  converges and $\lim_{t\rightarrow +\infty}T(r(t)) = \hat{T}$.
\end{Theorem}

\subsection{Open problems}
To conclude this section, we pose two natural questions:\\\\
\textbf{Question 1:} Can our current results be extended to the setting $\Theta: E \rightarrow$ $[0, \pi)$, which is known to be established in hyperbolic background geometry ?\\\\
\textbf{Question 2:} Does an analogous version of Theorem \ref{Thm He} hold in spherical background geometry ?\\\\

\section{Preliminaries}
\subsection{Circle patterns in spherical background geometry}
\begin{figure}[htbp]
	\centering          
	\includegraphics[width=0.8\textwidth]{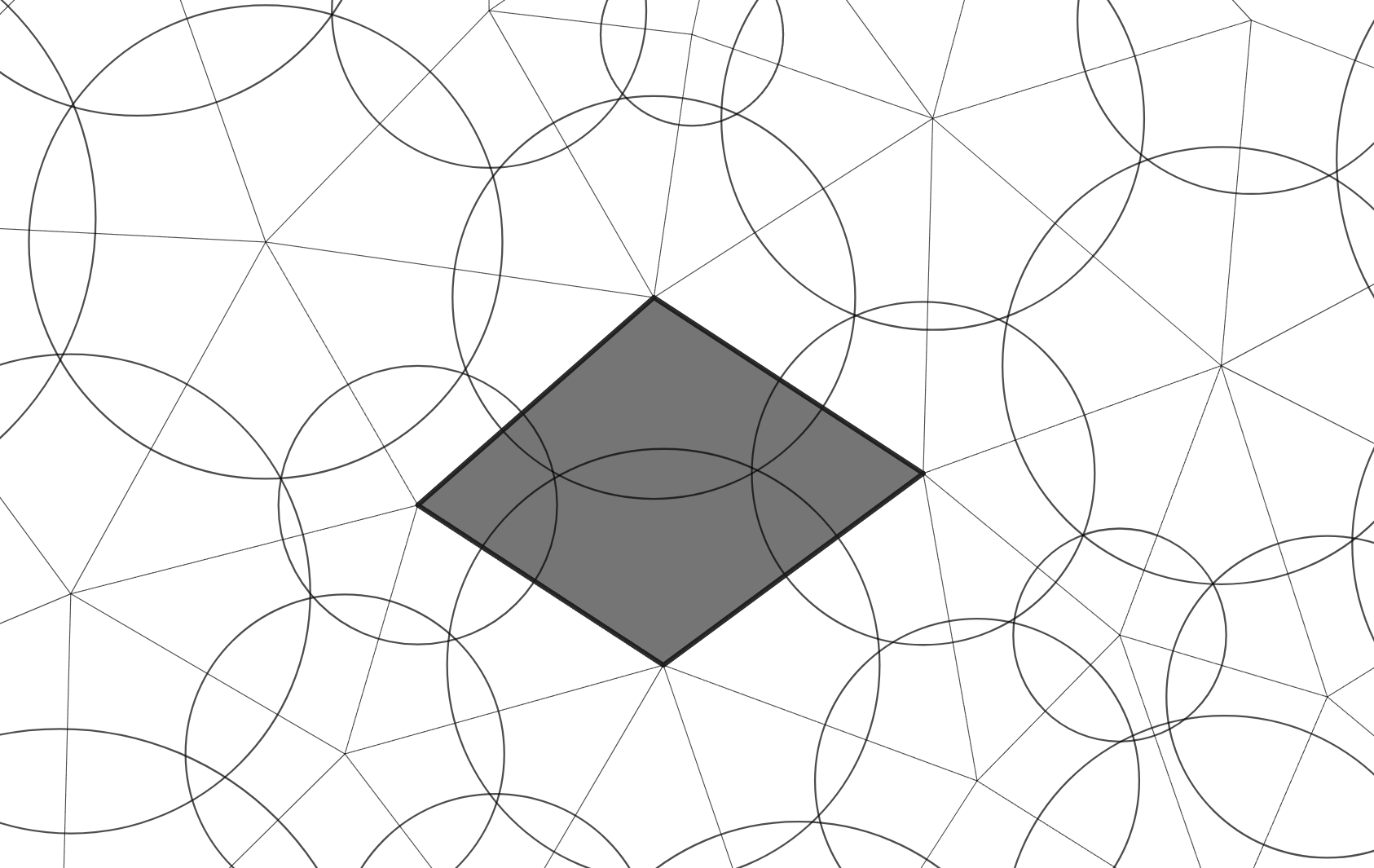} 
	\caption{A circle pattern on $S$}   
	\label{figure1}  
\end{figure}

Let $S$ be a surface equipped with an infinite cellular decomposition $\mathcal{D}=(V, E, F)$, where $V, E$ and $F$ represent the sets of vertices, edges and faces, respectively (see Figure \ref{figure1}). For a surface S endowed with a metric $\mu$, we write $(S, \mu)$. In this paper, we assume that $\mathcal{D}$ is locally finite, i.e. $\deg (v)<+\infty$ for any $v \in V$. For simplicity, we denote the vertex set $V$ as $\{v_i\}_{i=1}^\infty$. Besides, the adjacency relation between vertices $v_i$ and $v_j$ is denoted by $v_i \sim v_j$ or $i \sim j$ when they are connected by an edge in $E$.

A circle pattern metric refers to a positive function $r: V \rightarrow (0, \frac{\pi}{2})$ assigned to the vertices. Additionally, an intersection angle is defined as a function $\Theta: E \rightarrow(0, \frac{\pi}{2}]$ that assigns an angle to each edge. Sometimes, for simplicity, we denote $\Theta(e)$ by $\Theta_{ij}$, where $e = [v_i, v_j] \in E$.

Given a triple $(\mathcal{D}, r, \Theta)$, for each edge $e = [v_i, v_j] \in E$, we construct a corresponding spherical quadrilateral $Q_e$ (see Figure \ref{figure1}) satisfying  
\[
\angle v_i A v_j = \angle v_i B v_j = \pi - \Phi(e), \quad |v_i A| = |v_i B| = r_i, \quad |v_j A| = |v_j B| = r_j.
\]
By the spherical law of cosines, the quadrilateral $Q_e$ is uniquely determined up to isometry (see \cite[Lemma 7]{Nie}).  
Gluing all such quadrilaterals $Q_e$ along the cellular decomposition $\mathcal{D}$ yields a cone spherical metric $\mu(\mathcal{D}, r, \Theta)$ on the surface $S$, possibly with conical singularities. For details of the gluing procedure, see \cite[Chapter 3]{BBI22}.

We denote the angle $\angle A v_i B$ by $\Theta(e, v_i)$. Moreover, we use the notation $v < e$ (or $e < f$) to indicate that a vertex $v$ is incident to an edge $e$ (or an edge $e$ is incident to a face $f$). For a vertex $v_i \in V$, the cone angle at $v_i$ is defined by
\[
\alpha_i = \sum_{e: v_i < e} \Theta(e, v_i).
\]

\begin{figure}[htbp]
	\centering          
	\includegraphics[width=0.6\textwidth]{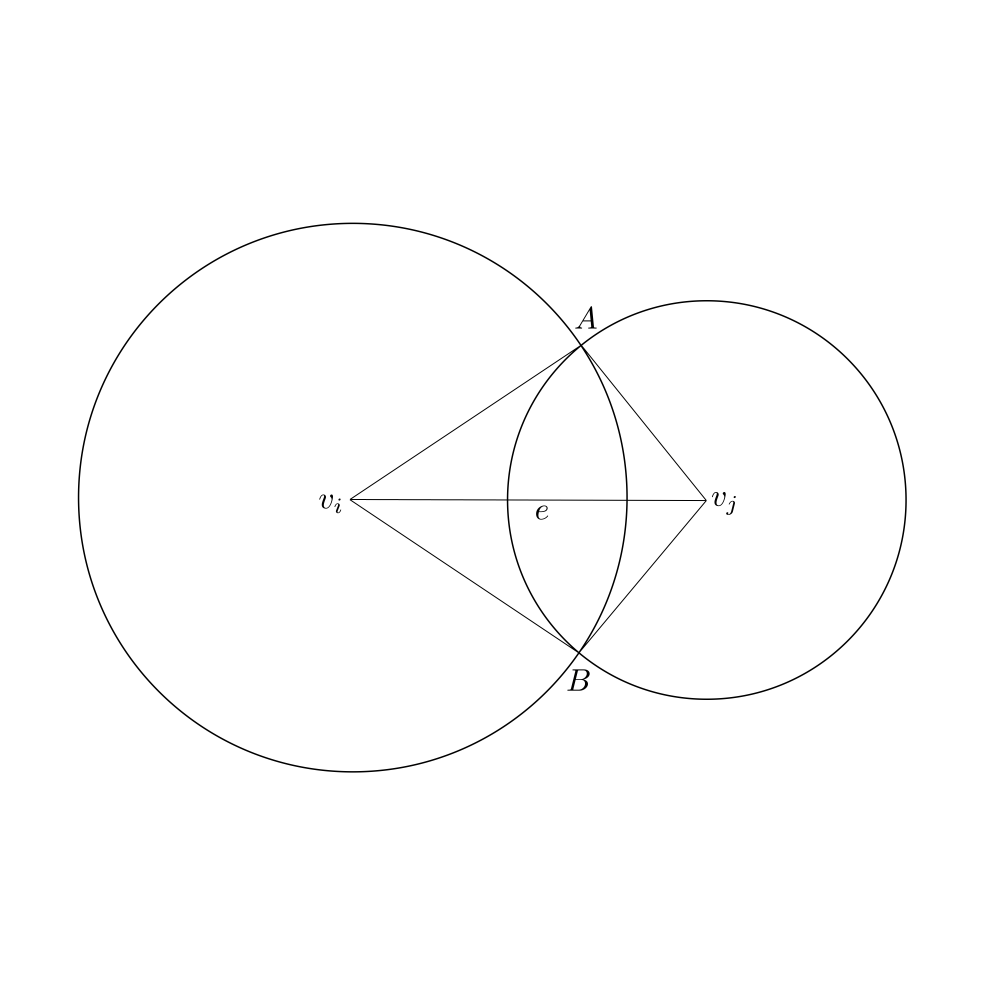} 
	\caption{Spherical quadrilateral $Q_e$}   
	\label{figure2}  
\end{figure}

Furthermore, let $C_i$ denote the circle centered at $v_i$ with radius $r_i$ on the surface $(S,\mu)$ equipped with the metric $\mu(\mathcal{D}, r, \Theta)$.
We denote by $k_i$ and $l_i$ the {geodesic curvature} and {circumference} of $C_i$, respectively. These quantities are given by
\[
k_i = \cot r_i, \quad l_i = \alpha_i \sin r_i.
\]
Integrating the geodesic curvature along $C_i$ yields the {total geodesic curvature} $T_i$ of $C_i$, which satisfies
\begin{equation}
	T_i = k_i l_i = \alpha_i \cos r_i.
\end{equation}

A basic problem in circle pattern theory is the prescribed total geodesic curvature problem. More precisely, given a set of intersection angles $\{\Theta(e)\}_{e \in E}$ and target total geodesic curvatures $\{\hat{T}_i\}_{v_i \in V}$, we want to find radii $\{r_i\}_{v_i \in V}$ ensuring that the geodesic curvature of each circle $C_i$ satisfies

$$
T_i=\hat{T}_i, \quad \forall v_i \in V .
$$
To solve this problem, Ge-Hua-Zhou \cite{GHZ24} introduce the following prescribed geodesic curvature flow:
\begin{equation}\label{flowequation}
	\frac{\mathrm{d} r_i}{\mathrm{~d} t}=\frac{(T_i-\hat{T}_i)}{2} \sin (2 r_i), \quad \forall v_i \in V .
\end{equation}

In \cite{GHZ24}, for a finite cellular decomposition $\mathcal{D}=(V, E, F)$, they prove that flow \eqref{flowequation} exists for all time. Moreover, it converges if and only if $\{\hat{T}_i\}_{v_i \in V}$ satisfy the condition
\begin{equation}
	\hat{T}_v>0, \quad \forall v \in V,
\end{equation}
and
\begin{equation}\label{1.4}
	\sum_{v \in U} \hat{T}_v<2 \sum_{e \in E(U)} \Theta(e), \quad \forall U \subset V.
\end{equation}

In our paper, we study the flow \eqref{flowequation} where $\mathcal{D}=(V, E, F)$ is an infinite cellular decomposition.

\subsection{ The potential function for total geodesic curvatures}
Let $u_i=\ln \cot r_i$ and $T_{(e, v)}=\Theta(e, v) \cos r_v$, then we have the following lemma, which can be found in Nie\cite{Nie}.
\begin{Lemma}\label{nie_lemma}
Given $(\mathcal{D}, r, \Theta)$ and $e=\{v_i, v_j\} \in E$, we have
	\begin{equation}
     \frac{\partial T_{(e, v_i)}}{\partial u_j}=\frac{\partial T_{(e, v_j)}}{\partial u_i},
	\end{equation}
and
\begin{equation}
	\frac{\partial T_{(e, v_i)}}{\partial u_j}<0,\quad \frac{\partial T_{(e, v_i)}}{\partial u_i}>0, \quad \frac{\partial(T_{(e, v_i)}+T_{(e, v_j)})}{\partial u_j}>0 .
\end{equation}
\end{Lemma}

\begin{Proof}
First, we have the cotangent 4-part formula in spherical background geometry, which is given by
\begin{equation}\label{L2.5}
	\cot \frac{\Theta(e, v_i)}{2}=\frac{1}{\sin \Theta(e)}(\cot r_j \sin r_i+\cos r_i \cos \Theta(e)).
\end{equation}
Then, differentiating both sides of \eqref{L2.5} with respect to $u_i$ and $u_j$, we have  
\begin{equation}
\frac{\partial T_{(e, v_i)}}{\partial u_j}=-\frac{2 \cos r_i \cos r_j \sin \frac{\Theta(e, v_i)}{2} \sin \frac{\Theta(e, v_j)}{2}}{\sin \Theta(e)}<0,
\end{equation}
and
\begin{equation}
\frac{\partial(T_{(e, v_i)}+T_{(e, v_j)})}{\partial u_i}=\sin ^2 r_j \cos r_i(\Theta(e, v_i)-\sin \Theta(e, v_i))>0 .
\end{equation}
where we use the law of sines, i.e.
\begin{equation}
	\frac{\sin \frac{\Theta(e, v_i)}{2}}{\sin r_j}=\frac{\sin \frac{\Theta(e, v_j)}{2}}{\sin r_i}.
\end{equation}
We omit these details; see \cite{GHZ24} for the complete proof.

\end{Proof}

Nie gives a geometric interpretation of Lemma \ref{nie_lemma} in \cite{Nie}. Given $e={v_1,v_2}\in E$, we define 
$$
\omega(e)=T_{(e, v_1)} \mathrm{d} u_1+T_{(e, v_2)} \mathrm{d} u_2 ,\quad (u_1,u_2) \in  \mathbb{R}^2,
$$
which is a closed form by Lemma \ref{nie_lemma}. Then, we can define 
$$
\mathcal{E}_{e}(u_1, u_2)=\int^{(u_1, u_2)} \omega(e),
$$
which is a convex potential function. Additionally, we can also define a potential function $\mathcal{E}_e$ for each edge $e$. Finally, we define the following potential function
$$
\mathcal{E}(u)=\sum_{e=\{v_i, v_j\} \in E} \mathcal{E}_e(u_i, u_j)-\sum_{v \in V} \hat{T}_v u_v ,
$$
for a given $\{\hat{T}_i\}_{v_i \in V}$. The potential function was introduced by Nie\cite{Nie} and has been widely used in many studies, such as \cite{GHZ24, LZ25}.

\subsection{A maximum principle for infinite graph}
Given an undirected infinite graph $G=(V, E)$, where $V$ and $E$ are the vertex and edge sets of $G$ respectively, we write $i \sim j$ if $v_i$ and $v_j$ are connected by an edge in $E$. For a function $f: V \rightarrow \mathbb R$, the  discrete Laplacian operator $\Delta_G$ is given by

\begin{equation}
\Delta_G f_i=\sum_{ j:j\sim i} \omega_{ij }(f_j-f_i)
\end{equation}
where $\omega_{ij}$ is the weight on $E$.

\begin{Lemma}\label{Maximum}
	(Maximum principle for infinite graph). Let $G=(V, E)$ be an undirected infinite graph, $\omega_{ij}(t)$ be weight on $E$ for $t \in[0, \tau]$ with $\tau>0$. And there is an uniform constant $C$ such that
	$$
	\sum_{j:j \sim i} \omega_{i j}(t) < C,
	$$
	for any $v_i \in V$ and for any $t \in[0, \tau]$. Given a function $f: V \times[0, \tau] \rightarrow \mathbb{R}$ satisfies
	$$
	\frac{df}{dt}\leq \Delta_{G} f+g f.
	$$
	If $f$ is a bounded function in $V \times[0, \tau]$ with $f(0) \leq 0$ and $g \leq C_0$ for some constant $C_0$, then
	$$f(t) \leq 0, \quad \forall t \in [0, \tau].$$
\end{Lemma}
\begin{Proof}
Let $\tilde{f}=e^{-C_0 f}$ and $g=\tilde{g}-C_0$, we have

\begin{equation}
\begin{aligned}
	\frac{d \tilde{f}}{d t} & =e^{-C_0 t} \frac{d f}{d t}-C_0 e^{-C_0 t} f \\
	& \leqslant e^{-C_0 t} \Delta_G f+e^{-C_0 t} g f-C_0 e^{-C_0 t} f \\
	& =\Delta_G \tilde{f}+\tilde{g} \tilde{f}.
\end{aligned}
\end{equation}
Hence, without loss of generality, we can assume $C_0=0$.

Fix a vertex $v_0 \in V$, we define the function $d_i^{[v_0]}=d^{[v_0]}(v_i)=\operatorname{dist}(v_i, v_0)$, where $\operatorname{dist}(v_i, v_0)$ is the distance between $v_0$ and $v_i$. Note that for any $v_i \in V$ and $t \in[0, \tau]$, we have

\begin{equation}\label{2.4}
\Delta_G d_i^{[v_0]} \leqslant \sum_{j: j \sim i} \omega_{i j}(t)<C.
\end{equation}
For $\delta>0$, we define  $ f^\delta=f-\delta d^{[v_0]}-C \delta t$ and by \eqref{2.4}, we have
\begin{equation}\label{2.5}
\begin{aligned}
	\frac{\mathrm{d} f^\delta}{\mathrm{d} t} & =  	\frac{\mathrm{d} f}{\mathrm{d} t}-C\delta \\
	&\leq \Delta_{G} f+g f-C\delta\\
	& =    \Delta_{G}  f^\delta + \delta\Delta_{G} d^{[v_0]}+ gf-C\delta \\
	& <\Delta_{G}  f^\delta + gf \\
	& \leq \Delta_{G}  f^\delta +  gf^\delta.
\end{aligned}
\end{equation}
Since $f$ is bounded, we have $f^\delta_i \rightarrow -\infty$ as $ d^{[v_0]}_i \rightarrow \infty$. Then we deduce that $f^\delta$ can attain its maximum on $V \times[0, \tau]$, and without loss of generality, we assume the maximum point is $(v_i, t_0)$. If $t_0= 0$, then $f_i^\delta(0) \leq f_i(0) \leq 0$; If $t_0> 0$, since $(v_i, t_0)$ is the maximum point, then we have $\frac{\mathrm{d} f^\delta_i}{\mathrm{d} t}( t_0) \geq 0$ and $\Delta_{G} f^\delta_i(t_0) \leq 0$. Hence, by \eqref{2.5} and $g\leq 0$, we obtain that $f^\delta_i(t_0)\leq 0$ which implies $f^\delta_i\leq 0$ on 
$V \times[0, \tau]$. Since $\delta$ is arbitrary, let  $\delta \rightarrow 0$, we obtain $f \leq 0$.

\end{Proof}

\section{Proof of main results}
Let $S$ be a surface equipped with a locally finite cellular decomposition $\mathcal D =(V, E,F)$. Then, we consider the following equation:
\begin{equation}\label{flow1}
	\frac{\mathrm{d} r_i}{\mathrm{~d} t}=\frac{(T_i-\hat{T}_i)}{2} \sin (2 r_i), \quad \forall i \in V .
\end{equation}
Let $u_i=\ln \cot r_i$, equation \eqref{flow1} is equivalent to 

\begin{equation}\label{flow2}
	\frac{\mathrm{d} u_i}{\mathrm{~d} t}=-(T_i-\hat{T}_i).
\end{equation}
Since $V$ is not finite, we could not use Picard-Lindelöf theorem to directly give the existence and uniqueness of equation \eqref{flow2}. Hence, our idea is to consider the solutions $u_n$ of equations \eqref{flow_finite} on finite subcomplexes of $\mathcal D$, and then use Arzelà-Ascoli theorem to prove that this sequence of solutions converges to the solution of equation \eqref{flow2}.

Specifically, we first choose a series of finite simple connected cellular decomposition, denoted by $\{\mathcal{D}^{[n]}=(V^{[n]}, E^{[n]}, F^{[n]})\}_{n=1}^{\infty}$, where

$$
\mathcal{D}^{[n]} \subset \mathcal{D}^{[n+1]}, \quad \bigcup_{n=1}^{\infty} \mathcal{D}^{[n]}=\mathcal{D}.
$$
Then, we consider the equations on $\mathcal{D}^{[n]}$ as follows
\begin{equation}\label{flow_finite}
	\begin{cases}\frac{d u_i^{[n]}(t)}{d t}=-(T_i^{[n]}-\hat{T}_i), \quad \forall v_i \in V^{[n]}, \forall t>0 \\ u_i^{[n]}(0)=u_i(0), \quad \forall v_i \in V^{[n]}\\ u_i^{[n]}(t)=u_i(0), \quad \forall v_i \notin V^{[n]}, \forall t \geq 0\end{cases}
\end{equation}
where $T_i^{[n]}(t)=T_i(u^{[n]}(t))$. By definition, for any $v_i \in V$, the vertex curvature $T_i^{[n]}$ only depends on $u_j^{[n]}$, where $v_j=v_i$ or $v_j \sim v_i$. Thus, the vertex curvature $T_i^{[n]}$ at any $v_i \in V^{[n]}$ depends only on $u_j^{[n]}$ for $v_j \in V^{[n]}$. Since $V^{[n]}$ is finite, by Picard-Lindelöf theorem, we derive the following existence and uniqueness theorem:

\begin{Lemma}\label{flow_finite_exi}
	The flow \eqref{flow_finite} has a unique solution $u^{[n]}(t)$ which exists for all time $t \geq 0$.
\end{Lemma}

\begin{Proof}
	Since each $T_i^{[n]}$ is a smooth function of $u^{[n]}$, it follows that $T_i^{[n]}$ is locally Lipschitz continuous. By the classical theory of ODE, there exists some $\epsilon>0$ such that equation \eqref{flow_finite} has a solution $u^{[n]}(t)$ on $[0, \epsilon)$. Additionally, the theory of ODE guarantees the uniqueness of the solution on its maximal interval of existence. Thus, it suffices to show that $u^{[n]}(t)$ remains bounded on any $[0, T)$ which implies the maximal existence interval of $u^{[n]}(t)$ is $[0, \infty)$.
	
	By Gauss-Bonnet theorem, we obtain that
	\begin{equation}\label{GB_T}
	A_{ij}=2 \Theta (e_{ij})-T_{ij}-T_{ji},
	\end{equation}
	where $A_{ij}$ is the area of the intersection part of two spherical disks, $e_{ij}= [i,j]$ and $T_{ij}$ is the total geodesic curvature of $v_i$ contributed by $v_j$. Moreover, we have
	\begin{equation}\label{3.5}
		T_i = \sum_{j\sim i} T_{ij}.
	\end{equation}
	
	By \eqref{GB_T}, the total geodesic curvature $T_i$ of the circle $C_i$ is less than $\sum_{i<e} 2 \Theta(e)$ for each vertex $i$. Then, we have $|T_i^{[n]}(t)-\hat{T}_i| \leq \sum_{i<e} 2 \Theta(e)+|\hat{T}_i|$ for any $t \in [0, T)$  and 
	\begin{equation}\label{3.6}
	\begin{aligned}
		u_i^{[n]}(t) & =u_i^{[n]}(0)+\int_0^t \frac{d u_i^{[n]}(s)}{d s} d s \\
		& \leqslant u_i^{[n]}(0)+\int_0^t(\sum_{i<e} 2 \Theta(e)+|\hat{T}_i|) d s \\
		& \leqslant u_i^{[n]}(0)+T(\sum_{i<e} 2 \Theta(e)+|\hat{T_i}|) \\
		& <+\infty .
	\end{aligned}
	\end{equation}
	Therefore, the flow \eqref{flow_finite} exists for $t \in[0, \infty)$.
\end{Proof}
Theorem \ref{exist} is equivalent to the following theorem:
\begin{Theorem}\label{existence}
	The flow\eqref{flow2}  exists for all time $t \geq 0$.
\end{Theorem}

\begin{Proof}
	Firstly, we prove that for any $\tau>0$ and any vertex $v_i \in V$, we have that $\{u_i^{[n]}(t)\}_{n=1}^{\infty}$ is a bounded sequence in $C^2[0, \tau]$. 
	
	Clearly, \( v_i \in V^{[n]}\)  for any $n\geq N$, where $N$ is a sufficiently large number. Denote $T_i^{[n]}(t) = T_i(u^{[n]}(t))$, then by \eqref{GB_T} and \eqref{3.5}, we have
	
	\begin{equation}\label{e1}
		|T_i^{[n]}(t)-\hat{T}_i| \leq \sum_{i<e} 2 \Theta(e)+|\hat{T}_i|, \quad t \in  [0, \tau] .
	\end{equation} 
	Furthermore, by \eqref{3.6}, we obtain that
	\begin{equation}\label{e2}
		|u_i^{[n]}(t)|\leq |u_i(0)|+  (\sum_{i<e} 2 \Theta(e)+|\hat{T}_i|)\tau , \quad t \in  [0, \tau]. 
	\end{equation}
	From \eqref{flow_finite}, combining \eqref{e1} and \eqref{e2}, we derive that 
	\begin{equation}\label{keyesti}
		\begin{aligned}
			\|u_i^{[n]}(t)\|_{C^1[0, \tau]} &= \sup _{t \in[0, \tau]}|u_i^{[n]}(t)|+\sup _{t \in[0, \tau]}|\frac{d u_i^{[n]}(t)}{d t}|\\
			&=\sup _{t \in[0, \tau]}|u_i^{[n]}(t)|+\sup _{t \in[0, \tau]}	|T_i^{[n]}(t)-\hat{T}_i| \\\
			&\leq |u_i(0)|+ (\sum_{i<e} 2 \Theta(e)+|\hat{T}_i|)(\tau+1).
		\end{aligned}
	\end{equation}
	Since each $T_i^{[n]}$ is a smooth function of $u^{[n]}$, $\frac{\partial T_i^{[n]}}{\partial u_i}$ and $\frac{\partial T_i^{[n]}}{\partial u_j}$ are also smooth functions of $u^{[n]}$. Hence, by \eqref{e2}, for $j \sim i$, we have that
	\begin{equation}\label{kij}
		\max \{|\frac{\partial T_i^{[n]}}{\partial u_i^{[n]}}|,|\frac{\partial T_i^{[n]}}{\partial u_j^{[n]}}|\}\leq C,
	\end{equation}
	where constant $C$ only depends on $\tau,i,j$ where $j \sim i$ and $ e\in E(v_i)$. Thus, by \eqref{e2} and \eqref{kij}, we have
	\begin{equation}\label{keyesti2}
		|\frac{d^2}{dt^2}u_i^{[n]}(t)|=|\frac{d}{dt} T_i^{[n]}(t)|=\sum_{j \sim i}|\frac{\partial T_i^{[n]}}{\partial u_j^{[n]}}|\cdot |\frac{d}{dt}u_j^{[n]}(t)| \leq  C'.
	\end{equation}
	Then, for any fixed vertex $v_i \in V$, we know that the sequence $\{u_i^{[n]}(t)\}_{n=1}^{\infty}$ is bounded in $C^2[0, \tau]$. Hence, by the Arzelà-Ascoli theorem, we know that there is a subsequence of $\{u_i^{[n]}(t)\}_{n=1}^{\infty}$, denote by $\{u_{i,\tau}^{[n^{(i,\tau)}_k]}(t)\}_{k=1}^{\infty}$  which converges to a function in $C^1[0, \tau]$. 
	Since $V$ is countable, without loss of generality, we assume $V=\mathbb{N}$. Besides, we denote $\{u_{i,j}^{[n]}(t)\}_{n=1}^{\infty}$=$\{u_{i,\tau}^{[n]}(t)\}_{n=1}^{\infty}$ where $\tau =j \in \mathbb{Z}$ and endow $\mathbb{N}^2$ with the diagonal order(see Figure \ref{figure3}):
	\begin{equation}
	(i, j) \prec(i^{\prime}, j^{\prime}) \quad \text { if and only if } \quad \begin{cases}i+j<i^{\prime}+j^{\prime}, & \text { or } \\ i+j=i^{\prime}+j^{\prime} & \text { and } i<i^{\prime}.\end{cases}
	\end{equation}
	
	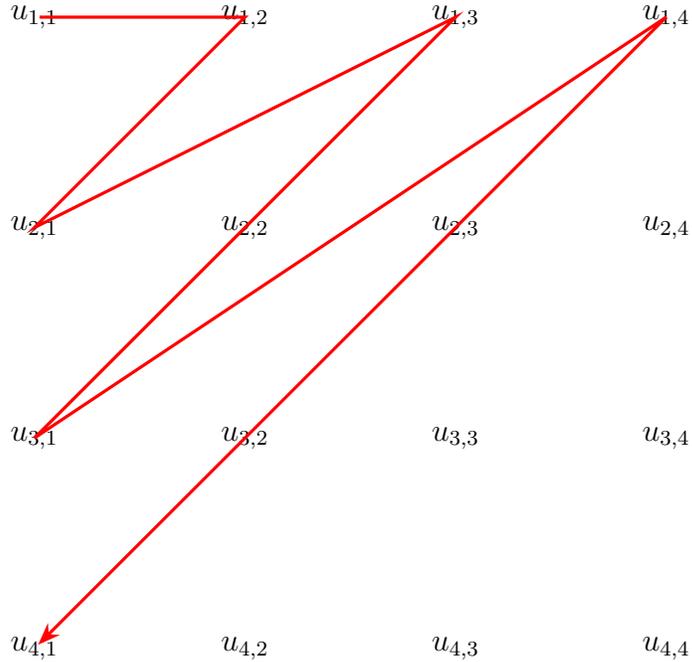
\begin{figure}[htbp]
		\centering 
		
		\begin{tikzpicture}[
			cell/.style={
				minimum size=2cm, 
				anchor=center,
				font=\large
			},
			arrow/.style={
				-Stealth, 
				red, 
				line width=1.2pt,
				shorten >=2pt, 
				shorten <=2pt
			}
			]
			\matrix (m) [
			matrix of math nodes, 
			nodes={cell}, 
			column sep=0.8cm,
			row sep=0.8cm
			] {
				u_{1,1} & u_{1,2} & u_{1,3} & u_{1,4} \\
				u_{2,1} & u_{2,2} & u_{2,3} & u_{2,4} \\
				u_{3,1} & u_{3,2} & u_{3,3} & u_{3,4} \\
				u_{4,1} & u_{4,2} & u_{4,3} & u_{4,4} \\
			};
			
			\draw[arrow] 
			(m-1-1.center) -- (m-1-2.center) -- (m-2-1.center)
			-- (m-1-3.center) -- (m-2-2.center) -- (m-3-1.center)
			-- (m-1-4.center) -- (m-2-3.center) -- (m-3-2.center) 
			-- (m-4-1.center);
		\end{tikzpicture}
		
		\caption{The diagonal order} 
		\label{figure3} 
	\end{figure}

	For any fixed $(i,\tau)$, by \eqref{keyesti} and \eqref{kij}, we know that $\{u_{i,\tau}^{[n]}(t)\}_{n=1}^{\infty}$ is bounded in $C^2[0, \tau]$.
	 Hence, we can inductively select such a subsequence: for $(i,j)=(1,1)$, we take $\{n^{(1,1)}_k\}_{k=1}^{\infty} \subseteq \mathbb{N}$ such that $\{u_{(1,1)}^{[n^{(1,1)}_k]}(t)\}_{k=1}^{\infty}$ is a convergent subsequence of $\{u_1^{[n]}(t)\}_{n=1}^{\infty}$ in $C^1[0, 1]$; 
	 for $(i,j)\succ (1,1)$, we take $\{n^{(i,j)}_k\}_{k=1}^{\infty} \subseteq \{n^{(i^\prime,j^\prime)}_k\}_{k=1}^{\infty}$ where if $i \geq 2$, then $(i, j)=(i-1, j+1)$; otherwise, $(i, j)=(j-1,1)$,  such that $\{u_{i,j}^{[n^{(i,j)}_k]}(t)\}_{k=1}^{\infty}$ is a convergent subsequence of $\{u_i^{[n^{(i^\prime,j^\prime)}_k]}(t)\}_{n=1}^{\infty}$ in $C^1[0, \tau]$. 
	 
	 Then, we take the diagonal sequence $\{n^{(j,j)}_j\}_{j=1}^{\infty}$. Fixed a vertex $v_i$ and $\tau =k$, note that $\{u_{i,k}^{[n^{(j,j)}_j]}(t)\}_{j=1}^{\infty}$ converges to a function $u_{i,j}(t)$ in $C^1[0, k]$. Hence, let $j \rightarrow+\infty$ in \eqref{flow_finite}, we deduce that $\{u_{i,k}(t)\}_{j=1}^{\infty}$ is the solution to equation \eqref{flow2} in $[0, k]$. Finally, we define
	 \begin{equation}\label{3.13}
	 u_i(t)=u_{i,k}(t), \quad t\leq k.
	 \end{equation}

     For any $k <k^\prime$, $u_{i,k}(t)$ is the limit of $\{u_{i,k}^{[n^{(i,k)}_j]}(t)\}_{j=1}^{\infty}$ in $C^1[0, k]$ and $u_{i,k^\prime}(t)$ is the limit of $\{u_{i,k^\prime}^{[n^{(i,k^\prime)}_j]}(t)\}_{j=1}^{\infty}$  in $C^1[0,k^\prime]$. Since $\{n^{(i,k^\prime)}_j\}_{j=1}^{\infty} \subseteq \{n^{(i,k)}_j\}_{j=1}^{\infty}$, we have $u_{i,k}(t)=u_{i,k^\prime}(t)$ on $[0,k]$, which implies the definition of \eqref{3.13} is well-defined. Since $u_{i,k}(t)$ is the solution to \eqref{flow2} on $[0,k]$ for any $k \in \mathbb{N}$, we deduce that $u_{i}(t)$ is the solution to \eqref{flow2} on $[0, \infty)$.
\end{Proof}

Theorem \ref{conv} is equivalent to the following theorem:
\begin{Theorem}
	For an initial value $u(0)$ such that $T(0)\geq \hat{T}$, then there exists a solution $u(t)$ that converges and $T(\infty) \rightarrow \hat{T}$.
\end{Theorem}
\begin{Proof}
	Let $u(t)$ be the flow which is obtained in Theorem \ref{existence} with the initial value $u(0)$ such that $T(0)\geq \hat{T}$.
	
	Firstly, we prove that $u(t)$ decreases.
	Recall that $T^{[n]}(t) = T(u^{[n]}(t))$, we define
	\begin{equation}
		f_i^{[n]}(t)= 
		\begin{cases}
			T_i^{[n]}(t) -\hat T_i & \text{if } i \in V^{[n]}, \\
			0 & \text{if }  i \notin V^{[n]}.
		\end{cases}
	\end{equation}
	If $v_i \in V^{[n]}$, since $u^{[n]}(t)$ is the solution to the flow \eqref{flow_finite}, we have
	\begin{equation}\label{pf1}
		\begin{aligned}
			\frac{d}{d t} f_i^{[n]}(t)=&-\frac{\partial f_i^{[n]}}{\partial u_i} \frac{d u_i^{[n]}(t)}{d t}-\sum_{j \sim i, v_j \in V^{[n]}} \frac{\partial f_i^{[n]}}{\partial u_j} \frac{d u_j^{[n]}(t)}{d t}(t) \\
			=&\frac{\partial f_i^{[n]}}{\partial u_i} f_i^{[n]}(t)-\sum_{j \sim i, v_j \in V^{[n]}} \frac{\partial f_i^{[n]}}{\partial u_j} f_j^{[n]}(t) \\
			=&-\sum_{j \sim i, v_j \in V^{[n]}} \frac{\partial f_i^{[n]}}{\partial u_j} (f_j^{[n]}(t)-f_i^{[n]}(t))-(\frac{\partial f_i^{[n]}}{\partial u_i} +\sum_{j \sim i, v_j \in V^{[n]}} \frac{\partial f_i^{[n]}}{\partial u_j} )f_i^{[n]}(t).
		\end{aligned}
	\end{equation}
	Otherwise, we have $\frac{d}{d t} f_i^{[n]}(t)=0$. Then, we define
	\begin{equation}
		\omega_{ij}(t)= 
		\begin{cases}
			-\frac{\partial f_i^{[n]}}{\partial u_j} & \text{if } i,j \in V^{[n]}, \\
			0 & \text{other cases,}
		\end{cases}
	\end{equation}
	and 
	\begin{equation}
		g(t) = 
		\begin{cases}
						-(\frac{\partial f_i^{[n]}}{\partial u_i} +\sum_{v_j \sim v_i, v_j \in V^{[n]}} \frac{\partial f_i^{[n]}}{\partial u_j}) & \text{if } i \in V^{[n]}, \\
			0 &\text{if } i \notin V^{[n]}.
		\end{cases}
	\end{equation}
	Hence, by $\frac{d}{d t} f_i^{[n]}(t)=0$ for $v_i \notin V^{[n]}$ and \eqref{pf1}, we have
	\begin{equation}
			\frac{df^{[n]}}{dt} = \Delta_{\omega(t)} f^{[n]}+g f^{[n]}.
	\end{equation}
	From Lemma \ref{nie_lemma}, we know $\omega_{ij}=\omega_{ji}\geq 0$. Given a $\tau>0$, for $t \in [0, \tau]$, since $\{\mathcal{D}^{[n]}=(V^{[n]}, E^{[n]}, F^{[n]})\}_{n=1}^{\infty}$ is  a finite cellular decomposition, by \eqref{keyesti2}, we know $|g|\leq C_0$ for an uniform constant $C_0$. Besides we also have $|\omega_{ij}|\leq C_1$ for an uniform constant $C_1$ by \eqref{keyesti2}. Thus, for any $v_i \in V$ and for any $t \in[0, \tau]$, there is an uniform constant $C$ such that
	\begin{equation}
	\sum_{j: j \sim i} \omega_{i j}(t) < C,
    \end{equation}
	Then, by Lemma \ref{Maximum}, we deduce that $f^{[n]}(t) \geq 0$, which implies $T_i(u^{[n]}(t))\geq \hat{T_i}$ for any $v_i \in V^{[n]}$. Since $\tau$ is arbitrary, we have $T(u(t))\geq \hat{T}$ for any $t \geq 0$. Then, by \eqref{flow2}, we deduce that $u(t)$ decreases.
	
	Secondly, we prove that $u_i(t)$ has the lower bound. If $u_i(t)$ is unbounded below, then there exists a sequence $\{t_k\}$ such that $u_i(t_k) \rightarrow -\infty$ as $k \rightarrow \infty$. Since $u_i=\ln \cot r_i$, we have $r_i(t_k) \rightarrow \frac{\pi}{2}$ as $k \rightarrow \infty$. Note that $T_i=\alpha_i \cos r_i$ and $\alpha_i \leq \pi \operatorname{deg}(v_i)$, then we have $T_i(t_k) \rightarrow 0$ as $k \rightarrow \infty$, which contradicts $T(u(t))\geq \hat{T}$. Hence, $u(t)$ converges.
	
	Finally, since $T(t)$ is a smooth function of $u(t)$, we have $T(t)$ converges. Then, by $T(u(t))\geq \hat{T}$, we know $T_i(\infty) \geq \hat{T_i}$. By contradiction, if $T_i(\infty) > \hat{T_i}$, then there exists a $\tau$ such that $T_i(t)\geq (T_i(\infty)+\hat{T_i})/2$ for any $t \geq \tau$. Then, by \eqref{flow2}, we have 
	\begin{equation}
    du_i(t)/dt \leq - (T_i(\infty)-\hat{T_i})/2<0
    \end{equation}
    when $t \geq \tau$. Since $ (T_i(\infty)-\hat{T_i})/2$ is a constant, we have $u_i(t) \rightarrow -\infty$ as $t \rightarrow \infty$, which contradicts with $u_i(t)$ has the lower bound.

\end{Proof}

\section{Acknowledgments}
Chang Li is supported by National Natural Science Foundation of China (Grant No. 12301079), the Fundamental Research Funds for the Central Universities and the Research Funds of Renmin University of China (Grant No. 24XNKJ01). Yangxiang Lu is supported by National Natural Science Foundation of China (Grant No. 12341102, No.12122119).
Hao Yu is supported by NSF of China (No. 12341102).

\noindent Chang Li, 
chang\_li@ruc.edu.cn\\
\emph{School of Mathematics, Renmin University of China, Beijing 100872, P. R. China}\\

\noindent Yangxiang Lu, 
2023000744@ruc.edu.cn\\
\emph{School of Mathematics, Renmin University of China, Beijing 100872, P. R. China}\\

\noindent Hao Yu, b453@cnu.edu.cn\\
\emph{Academy for multidisciplinary studies, Capital Normal University, Beijing, 100048, P. R. China}


\begin{thebibliography}{9}
	
	
	
	
	
	\bibitem{BBI22}
	D. Burago, Y. Burago, S. Ivanov. A course in metric geometry, volume 33. American Mathematical Society, 2022.
	

	                
	                
	\bibitem{CL03}  
	B. Chow, F. Luo. Combinatorial Ricci flows on surfaces. J. Differential Geom, 2003, 63(1): 97-129.
	
	\bibitem{CJN21}
	J. Cheeger, W. Jiang, A. Naber. Rectifiability of singular sets of noncollapsed limit spaces with Ricci curvature bounded below. Ann. Math. 193 (2021), 407-538.
	
	\bibitem{Ve91}
	Y. Colin de Verdi\`{e}re. Un principe variationnel pour les empilements de cercles. Invent. Math. 1991, 104(1): 655-669.
	
	\bibitem{FGH21}
	K. Feng, H. Ge, B. Hua. Combinatorial Ricci flows and the hyperbolization of a class of compact 3-manifolds. Geom. Topol. 26 (2022) 1349-1384.
	

	
	\bibitem{Ge12}
	H. Ge. Combinatorial methods and geometric equations. Thesis (Ph.D.)-Peking University, 2012.
	
	\bibitem{Ge17}
	H. Ge. Combinatorial Calabi flows on surfaces. Transactions of the American Mathematical Society, 370(2):1377-1391, 2017.
	
	\bibitem{GH18}
	H. Ge and B. Hua. On combinatorial Calabi flow with hyperbolic circle patterns.
	Advances in Mathematics, 333:523-538, 2018.
	
	\bibitem{GH20}
	H. Ge, B. Hua. 3-dimensional combinatorial Yamabe flow in hyperbolic background geometry. Trans. Amer. Math. Soc. 2020, 373(7): 5111-5140.
	
	
	\bibitem{GHZ21}
	H. Ge, B. Hua, Z. Zhou. Combinatorial Ricci flows for ideal circle patterns. Adv. Math. 2021, 383: 107698.
	
	\bibitem{GHZ21-2}
	H. Ge, B. Hua, Z. Zhou. Circle patterns on surfaces of finite topological type. Amer. J. Math. 2021, 143(5): 1397-1430.
	
	\bibitem{GHZ24}
	H. Ge, B. Hua, P. Zhou. A combinatorial curvature flow in spherical background geometry. J. Funct. Anal. 286 (2024), no. 7, Paper No. 110335, 12 pp.
	
	\bibitem{GHZ25}
	H. Ge, B. Hua, P. Zhou. Combinatorial Ricci flows on infinite disk triangulations.
	arXiv:2504.05817v1.
	
	\bibitem{GJ16}
	H. Ge, W. Jiang. On the deformation of discrete conformal factors on surfaces. Calc. Var. Partial Differential Equations, 2016, 55(6): 1-14.
	
	
	\bibitem{GJ17}
	H. Ge, W. Jiang. On the deformation of inversive distance circle patterns, II. J. Funct. Anal. 2017, 272(9): 3573-3595.
	
	\bibitem{GJ17-2}
	H. Ge, W. Jiang. On the deformation of inversive distance circle patterns, III. J. Funct. Anal. 2017, 272(9): 3596-3609.
	
	\bibitem{GJS22}
	H. Ge, W. Jiang, L. Shen. On the deformation of ball patterns. Adv. Math. 2022, 398: 108192.
	
	\bibitem{GX16DGA}
	H. Ge, X. Xu. 2-dimensional combinatorial Calabi flow in hyperbolic background geometry. Differential Geometry and its Applications, 47:86-98, 2016.
	
	\bibitem{GX16}
	H. Ge, X. Xu. $\alpha$-curvatures and $\alpha$-flows on low dimensional triangulated manifolds. Calc. Var. Partial Differential Equations, 2016, 55(1).
	
	\bibitem{GX18}
	H. Ge, X. Xu. On a combinatorial curvature for surfaces with inversive distance circle pattern metrics. J. Funct. Anal. 2018, 275(3): 523-558.

	
	\bibitem{H82}
	R. Hamilton. Three-manifolds with positive Ricci curvature. J. Differential Geom. 1982, 17(2): 255-306.
	
	\bibitem{He95}
	Z. He, O. Schramm. Hyperbolic and parabolic patterns. Discrete Comput. Geom., 14(2):123-149, 1995.
	
	\bibitem{He99}
	Z. He. Rigidity of infinite disk patterns. Ann. of Math. (2), 149(1):1-33, 1999.
	
	\bibitem{JN21}
    W. Jiang, A. Naber. $L^2$ curvature bounds on manifolds with bounded Ricci curvature. Ann. Math. 193 (2021), 107-222.
	
	\bibitem{LZ25}
	Z. Lei, P. Zhou. Combinatorial Calabi flows for ideal circle patterns in spherical background geometry. J. Differential Equations 427 (2025), 676-688.
	
	\bibitem{MR89}
	A. Marden, B. Rodin. On Thurston's formulation and proof of Andreev's theorem. Computational Methods and Function Theory, LNM 1435, Springer-Verlag, New York, $1989,103-115$.
	
	\bibitem{Nie}
	X. Nie. On circle patterns and spherical conical metrics. Proc. Amer. Math. Soc., 152(2024), 843-853.
	
	\bibitem{Perelman}
	G. Perelman. The entropy formula for the Ricci flow and its geometric applications. arXiv:math.DG/0211159v1.
	
	\bibitem{ro2}
	B. Rodin,  D. Sullivan. The convergence of circle pattern to the Riemann mapping. J Differential Geom, 26, 1987, 349-360.

	\bibitem{st}
	K. Stephenson.  Introduction to circle pattern. The theory of discrete analytic functions, Cambridge University Press, Cambridge, 2005.
	
	\bibitem{Th76}
	W. Thurston. Geometry and topology of three-manifolds. Princeton Lecture Notes, 1976.
	
	\bibitem{Th80}
	W. Thurston. The geometry and topology of 3-manifolds, Chapter 13, Princeton Univ. Math. Dept. Notes. Princeton, New Jersey, 1980.
	
\end{thebibliography}
\end{document}